\documentclass[12pt]{article}

\usepackage{amsmath}
\usepackage{amsfonts}
\usepackage{amsthm}

\begin{document}

\newcommand{\A}{\mbox{${{{\cal A}}}$}}


\author{Attila Losonczi}
\title{Means of infinite sets I}

\date{6 June 2018}

\newtheorem{thm}{\qquad Theorem}[section]
\newtheorem{prp}[thm]{\qquad Proposition}
\newtheorem{lem}[thm]{\qquad Lemma}
\newtheorem{cor}[thm]{\qquad Corollary}
\newtheorem{rem}[thm]{\qquad Remark}
\newtheorem{ex}[thm]{\qquad Example}
\newtheorem{df}[thm]{\qquad Definition}
\newtheorem{prb}{\qquad Problem}

\maketitle

\begin{abstract}

\noindent

We open a new field on how one can define means on infinite sets. We investigate many different ways on how such means can be constructed. One method is based on sequences of ideals, other deals with accumulation points, one uses isolated points, other deals with average using integral, other with limit of average on surroundings and one deals with evenly distributed samples. We study various properties of such means and their relations to each other.

\noindent
\footnotetext{\noindent
AMS (2010) Subject Classifications: 26E60, 28A10, 28A78  \\

}

\end{abstract}

\section{Introduction}
In this paper we are going to study the ways of how can one generalize the arithmetic mean for an infinite bounded subset of $\mathbb{R}$. As well known one can calculate the arithmetic mean for finite sets and there is a straightforward generalization for sets with finite positive Lebesgue measure (see Def \ref{davg}). One may ask if we can extend these methods in between or more generally in what way one can define (natural) means on infinite subsets of $\mathbb{R}$.

In this paper our aim is to find reasonably good and natural means for infinite bounded sets. Then study their properties and relations among them. We are going to present many methods where in some of them we deal with countable sets only.

Most of the methods described here can be easily generalized to quasi-arithmetic means as well or to more general means, however we are not going to deal such generalizations now. In this paper we focus on arithmetic type means only.

\medskip

We are planning a second paper on this topic. While this current paper mainly deals with constructing means and investigate their properties, the second paper is going to focus mainly on building and analysing concepts of this new field. 

\subsection{Basic notions and notations}
Throughout this paper function $\A()$ will denote the arithmetic mean of any number of variables.

\begin{df}\label{davg}If $H$ Lebesgue measurable, $\lambda(H)>0$ then $$Avg(H)=\frac{\int\limits_H x\ d\lambda}{\lambda(H)}.$$
\end{df}

For $K\subset\mathbb{R},\ y\in\mathbb{R}$ let us use the notation $$K^{y-}=K\cap(-\infty,y],\ K^{y+}=K\cap[y,+\infty).$$

Let $T_s$ denote the reflection to point $s\in\mathbb{R}$ that is $T_s(x)=2s-x\ (x\in\mathbb{R})$.

If $H\subset\mathbb{R},x\in\mathbb{R}$ then set $H+x=\{h+x:h\in H\}$. Similarly $\alpha H=\{\alpha h:h\in H\}\ (\alpha\in\mathbb{R})$.

$cl(H), H'$ will denote the closure and accumulation points of $H\subset\mathbb{R}$ respectively. Let $\varliminf H=\inf H',\ \varlimsup H=\sup H'$ for infinite bounded $H$.

\begin{df}A \textbf{generalized mean} is a function ${\cal{K}}:C\to \mathbb{R}$ where $C\subset P(\mathbb{R})$ consists of some (finite or infinite) bounded subsets of $\mathbb{R}$ and $\inf H\leq {\cal{K}}(H)\leq\sup H$ holds for all $H\in C$. We call ${\cal{K}}$  an \textbf{ordinary mean} if $C$ consists of finite sets only.
\end{df}

In the definition the required condition is an obvious generalization of internality on finite sets.

In this paper when we use the term "mean" we always refer to a generalized mean.
Usually ${\cal{K}},{\cal{M}}$ will denote means and $Dom({\cal{K}})$ denotes the domain of ${\cal{K}}$.

\subsection{Brief summary of the main notions and results}

Before we start to construct means we investigate what properties a generalized mean can have. In this current paper we just define some basic properties while the second paper deals with a much broader range of attributes. Here we have internality and strong internality that is the value of a mean has to be between the infimum and supremum of the set that is a defining relation for being a mean. We may also expect that a mean has to be invariant under some geometric transformations such that translation, reflection and contraction/dilation. We deal with one type of monotonicity and finite-independence that is the mean is invariant under adding or moving finite sets from the set. Some special properties are convexity and how a mean is invariant under closure or taking accumulation points.

In the paper we define 6 generalized means (however one is a group of means, not a single mean). Here we describe them roughly.

A mean can be defined by using the isolated points of a set. If the isolated points approximate the set in a sense that their arithmetic mean may approximate a mean  - that is the idea behind. This mean is strongly internal, monotone, convex and closed but not accumulated.

If the set of accumulation points is finite then its arithmetic mean may function as a mean of the set. If not, then we can take the accumulation points of the the accumulation points and so on. If any of these is finite then use that to define a mean. This mean is strongly internal, monotone, convex, closed and accumulated.

We can use ideals to define a kind of end points of sets, points over which the remains of the set is in the ideal i.e. small, in other words besides the end points the set is small. The arithmetic mean of the end points may function as a mean. This mean is monotone and convex.

There is a natural generalization of Avg using Hausdorff dimension and measure which somehow lay a bridge between \A\ and the Avg defined in \ref{davg}. This mean is strongly internal, strong monotone, convex but not closed, not accumulated.

If we have a set with 0 Lebesgue measure then Avg (\ref{davg}) cannot be applied. However we can take its $\epsilon$ surroundings for which Avg can be used and then take smaller and smaller $\epsilon$-s and see if the limit exists. If yes, then consider it as a mean. This mean is strongly internal and closed but not accumulated.

Finally we can take finite sample points from the set and calculate their arithmetic mean and we consider it as an approximation for a mean. Then take better and better samples and see if the limit exits. If yes, we found a new mean. It is important that the samples have to be evenly distributed. This mean is strongly internal, monotone and convex.

\section{Basic attributes of generalized means}
\subsection{Expected properties of a generalized mean}\label{ssexpprop}
Throughout these subsections ${\cal{K}}$ will denote a generalized mean. 

\smallskip

Usually we expect $Dom({\cal{K}})$ to be closed under finite union and intersection. Moreover closed under translation, reflection and contraction/dilation.

\smallskip

Most of our means $\cal{K}$ will be the extension of \A\ that is for finite sets it gives the arithmetic mean of the elements. Nevertheless we allow a mean to have a domain consists of infinite sets only.

\begin{df}\hfill

\begin{itemize}
\item $\cal{K}$ is called \textbf{internal} if for all $H\in Dom({\cal{K}})\ \inf H\leq {\cal{K}}(H)\leq\sup H$. 
${\cal{K}}$ is \textbf{strong internal} if for all infinite $H\in Dom({\cal{K}})$ $$\varliminf H\leq {\cal{K}}(H)\leq\varlimsup H.$$ 

\item $\cal{K}$ is \textbf{monotone} if $\sup H_1\leq\inf H_2$ implies that ${\cal{K}}(H_1)\leq {\cal{K}}(H_1\cup H_2)\leq {\cal{K}}(H_2)$. $\cal{K}$ is \textbf{strong monotone} if $\cal{K}$ is strong internal and $\varlimsup H_1\leq\varliminf H_2$ implies that ${\cal{K}}(H_1)\leq {\cal{K}}(H_1\cup H_2)\leq {\cal{K}}(H_2)$. 

\item The mean is \textbf{shift invariant} if $x\in\mathbb{R}, H\in Dom({\cal{K}})$ then $H+x\in Dom({\cal{K}}),\ {\cal{K}}(H+x)={\cal{K}}(H)+x$.

\item ${\cal{K}}$ is \textbf{symmetric} if $H\in Dom({\cal{K}})$ bounded and symmetric ($\exists s\in\mathbb{R}\ \forall x\ s+x\in H\Leftrightarrow s-x\in H$) implies ${\cal{K}}(H)=s$.

\item $\cal{K}$ is \textbf{homogeneous} if $H\in Dom({\cal{K}})$ then $\alpha H\in Dom({\cal{K}}),\ {\cal{K}}(\alpha H)=\alpha {\cal{K}}(H)$.

\item $\cal{K}$ is \textbf{finite-independent} if $H\in Dom({\cal{K}})$ is infinite, $V$ is finite then $H\cup V,H-V\in Dom({\cal{K}})$ and ${\cal{K}}(H)={\cal{K}}(H\cup V)={\cal{K}}(H-V)$.
\end{itemize}
\end{df}

A simple example of a mean ${\cal{M}}^{lis}$ on all bounded sets is 
$${\cal{M}}^{lis}(H)=\begin{cases}
\A(H) & \text{ if } H \text{ is finite }\\
\frac{\varliminf H+\varlimsup H}{2} & \text{ otherwise.}
\end{cases}$$
In subsection \ref{ssmbi} we will verify that it has all mentioned expected properties.


\begin{prp}If ${\cal{K}}$ is strong internal and $H'=\{h\}$ then ${\cal{K}}(H)=h$.
\end{prp}
\begin{proof}$h=\varliminf H\leq {\cal{K}}(H)\leq\varlimsup H=h$.
\end{proof}

\begin{prp}\label{pifiisi}If ${\cal{K}}$ is internal, finite-independent then $\cal{K}$ is strongly internal. 
\end{prp}
\begin{proof}Let $H\subset\mathbb{R},\ \epsilon>0$. Then ${\cal{K}}(H)={\cal{K}}(H\cap (-\infty,\varlimsup H+\epsilon])\leq \varlimsup H+\epsilon$ because we left out finitely many points. Since $\epsilon$ was arbitrary we get that ${\cal{K}}(H)\leq\varlimsup H$. Similar argument can be applied to $\varliminf$.
\end{proof}

\subsection{Some other properties}
\begin{df}\hfill
\begin{itemize}
\item ${\cal{K}}$ is \textbf{convex} if $I$ is a closed interval and ${\cal{K}}(H)\in I,\ L\subset I,H\cup L\in Dom({\cal{K}})$ then ${\cal{K}}(H\cup L)\in I$. 

\item ${\cal{K}}$ is called \textbf{closed} if $H,cl(H)\in Dom({\cal{K}})$ then ${\cal{K}}(cl(H))={\cal{K}}(H)$.

\item ${\cal{K}}$ is called \textbf{accumulated} if $H,H'\in Dom({\cal{K}})$ then ${\cal{K}}(H')={\cal{K}}(H)$. 
\end{itemize}
\end{df}

Obviously property "accumulated" is equivalent with that ${\cal{K}}(H)={\cal{K}}(H')={\cal{K}}(H'')={\cal{K}}(H''')=\dots$ if all sets are in $Dom({\cal{K}})$.

\medskip

We will often use the following simple fact.

\begin{lem}\label{laconvex}\A\ is convex.
\end{lem}
\begin{proof}Let $I$ be a closed interval, $H,L$ are finite and $\A(H)\in I,\ L\subset I$. It is known that if $A,B$ are disjoint finite sets with cardinality $|A|=a,|B|=b$
$$\A(A\cup B)=\frac{a\A(A)+b\A(B)}{a+b}=\frac{a}{a+b}\A(A)+\frac{b}{a+b}\A(B)$$
that is the convex combination of $\A(A)$ and $\A(B)$ hence it is between $\A(A)$ and $\A(B)$.

Now apply this to $H$ and $L-H$. For both $\A(H),\A(L-H)\in I$ hence so is $\A(H\cup L)$.
\end{proof}


We could consider a mean as Hausdorff continuous if it was a continuous function according to the Hausdorff pseudo-metric. However it cannot happen.

\begin{ex}Let $\cal{C}$ denote the set of all compact sets on $\mathbb{R}$ equipped with the Hausdorff metric, ${\cal{C}}_0\subset\cal{C}$ be the finite sets. Let $C=\{0,1\}$,\ $C_n=\{\frac{1}{n};1+\frac{1}{n};1+\frac{1}{2n}\}$. Clearly $C_n\to C$ in the Hausdorff metric but $\A(C_n)\to\frac{2}{3},\ \A(C)=\frac{1}{2}$. Hence \A\ is not a continuous function from ${\cal{C}}_0$ to $\mathbb{R}$.
\end{ex}

If ${\cal{K}}$ is not defined on finite sets then we still cannot expect ${\cal{K}}$ to be continuous assuming that it is strongly internal. 

\begin{ex}Let $H_1=\{0,\frac{1}{n}:n\in\mathbb{N}\}\cup\{1,1+\frac{1}{n}:n\in\mathbb{N}\}$.  Obviously ${\cal{K}}(H_1)\in[0,1]$. Let $L_{2k}=\{\frac{1}{n}:n\leq k\}\cup\{1,1+\frac{1}{n}:n\in\mathbb{N}\},L_{2k+1}=\{0,\frac{1}{n}:n\in\mathbb{N}\}\cup\{1+\frac{1}{n}:n\leq k\}$. Then $L_k\to H_1$ and $\forall k\ {\cal{K}}(L_{2k})=1,{\cal{K}}(L_{2k+1})=0$.
\end{ex}

\section{Simple generalized means}
\subsection{Mean by isolated points}
If the isolated points determine the set in the sense that $cl(H-H')=H$ then a mean can be defined by them using that for $\forall\delta>0\ H-S(H',\delta)$ is finite.

\begin{df}If $cl(H-H')=H$ then let $${\cal{M}}^{iso}(H)=\lim_{\delta\to 0+0}\A(H-S(H',\delta))$$ if it exists. 
\end{df}

\begin{lem}\label{ledsm}Let $(H_n),(L_n)$ be two infinite sequences of finite sets such that all sets are uniformly bounded, $\forall n\ H_n\cap L_n=\emptyset$ and $\A(H_n)\to a$. Moreover $\lim_{n\to\infty}\frac{|L_n|}{|H_n|}=0$. Then $\A(H_n\cup L_n)\to a$.
\end{lem}
\begin{proof}Clearly $$\A(H_n\cup L_n)=\frac{\sum\limits_{h_i\in H_n}h_i+\sum\limits_{h_j\in L_n}h_j}{|H_n\cup L_n|}=\frac{|H_n|}{|H_n\cup L_n|}\A(H_n)+\frac{|L_n|}{|H_n\cup L_n|}\A(L_n).$$ $\A(L_n)$ is bounded, $\frac{|L_n|}{|H_n\cup L_n|}\to 0$ and $\frac{|H_n|}{|H_n\cup L_n|}\to 1$ give the statement.
\end{proof}

\begin{thm}${\cal{M}}^{iso}$ is a generalized mean. Moreover it is finite-independent, strongly internal, monotone, shift invariant, symmetric, homogeneous, convex and closed. 
\end{thm}
\begin{proof}Clearly ${\cal{M}}^{iso}$ is internal since $H-S(H',\delta)\subset[\inf H,\sup H]$. 

It is also finite-independent because $H-H'$ is infinite and removing (or adding) finitely many new points would not change the limit. In order to prove that let $H_n=H-S(H',\delta),\ L_n=\{$the new points in $H-S(H',\delta)\}$. Then apply \ref{ledsm}.

Strong internality then follows from \ref{pifiisi}.

Let us show monotonicity. If $\sup H_1\leq\inf H_2$ then $$H_1\cup H_2-S((H_1\cup H_2)',\delta)=(H_1-S(H_1',\delta))\cup (H_2-S(H_2',\delta))$$
which gives that $\A(H_1-S(H_1',\delta))\leq\A(H_1\cup H_2-S((H_1\cup H_2)',\delta))$. When taking the limit we end up with ${\cal{M}}^{iso}(H_1)\leq{\cal{M}}^{iso}(H_1\cup H_2)$. The other inequality is similar.

To prove that ${\cal{M}}^{iso}$ is shift invariant, symmetric, homogeneous, it is enough to refer to the fact that $H-S(H',\delta)$ and \A\ both have the same properties.

To verify convexity let $I$ be a closed interval, ${\cal{M}}^{iso}(H)\in I,\ L\subset I,\ L,L\cup H\in Dom\ {\cal{M}}^{iso}$. It is known that if $A,B$ are disjoint finite sets with cardinality $|A|=a,|B|=b$
$$\A(A\cup B)=\frac{a\A(A)+b\A(B)}{a+b}$$
that is the convex combination of $\A(A)$ and $\A(B)$.
If we apply it for $A_{\delta}=H-S(H',\delta),\ B_{\delta}=\big(H\cup L-S((H\cup L)',\delta)\big)-A_{\delta}\subset L-S(L',\delta)$ then $\A(A_{\delta})\to p\in I,\ \A(B_{\delta})\in I$ hence in limit $(\delta\to 0+0)$ we get that $\A(A_{\delta}\cup B_{\delta})\to q\in I$ using that the limit exists because $H\cup L\in Dom\ {\cal{M}}^{iso}$.

To show that ${\cal{M}}^{iso}$ is closed it is enough to mention that $H$ and $cl(H)$ have the same set of isolated points.
\end{proof}

\begin{ex}For $H=\{0,1\}\cup\{\frac{1}{n}:n\in\mathbb{N}\}\cup\{1+\frac{1}{2^n}:n\in\mathbb{N}\}$, ${\cal{M}}^{iso}(H)=0$. 
\end{ex}
\begin{proof}Evidently $H'=\{0,1\}$. If $\delta=\frac{1}{k}$ then $$H-S(H',\delta)=\{\frac{1}{n}:n<k\} \cup \{1+\frac{1}{2^n}:2^n<k\}.$$
If we apply \ref{ledsm} for $H_k=\{\frac{1}{n}:n<k\},\ L_k=\{1+\frac{1}{2^n}:2^n<k\}$ then we get the statement.
\end{proof}

\begin{ex}${\cal{M}}^{iso}(H)$ does not exist always. 
\end{ex}
\begin{proof}Define a set in the following way. Let $H_1=\{1.7\}$. If $H_1,\dots,H_{n-1}$ are already defined then let $H_n$ consists of some finitely many points such that
$$H_n\subset
\begin{cases}
\big(\frac{1}{n+1},\frac{1}{n}\big)&\text{if n is even}\\
1+\big(\frac{1}{n+1},\frac{1}{n}\big)&\text{if n is odd}
\end{cases}$$ 
and $\A(H_1\cup\dots\cup H_n)\leq\frac{1}{4}$ when $n$ is even, $\A(H_1\cup\dots\cup H_n)\geq\frac{3}{4}$ when $n$ is odd. Then let $H=\bigcup\limits_{i=1}^{\infty}H_i$.

We then ended up with an infinite set $H\subset[0,2]$ such that $H'=\{0,1\}$ and $\A(H-S(H',\delta))$ can be smaller than $\frac{1}{4}$ or greater than $\frac{3}{4}$ depending on $\delta$ hence the limit does not exists.
\end{proof}

\begin{thm}${\cal{M}}^{iso}(H)$ is not accumulated.
\end{thm}
\begin{proof}It is easy to construct a set $H\subset[0,1]$ such that $H'=\{0\}\cup\{\frac{1}{n}:n\in\mathbb{N}\}$ and  $\forall\delta\ \A(\{h\in H:h>\delta\})\geq 0.5$. For that set we get ${\cal{M}}^{iso}(H)\geq 0.5,\ {\cal{M}}^{iso}(H')=0$.

For constructing such set let $H_1=\{1.5\}$. If $H_1,\dots,H_{n-1}$ are already defined then let $H_n$ consists of some points such that $H_n\subset S\big(\{\frac{1}{k}:k<n\},\frac{1}{n}\big),\forall k<n\ H_n\cap S(\frac{1}{k},\frac{1}{n})\ne\emptyset$ and $\A(H_1\cup\dots\cup H_n)\geq 0.5$. Obviously it can be done since we can add as many points around $1$ as we want. Then let $H=\bigcup\limits_{i=1}^{\infty}H_i$.
\end{proof}

\subsection{Mean by accumulation points}

Let us recall the classic definition. $H^{(0)}=H,\ H^{(1)}=H'$ where $H'$ denotes the accumulation points of H. Then $H^{(n+1)}=(H^{(n)})'\ (n\geq 0)$.

\medskip

Assume that $H$ is infinite bounded. Then there are two cases. Either there is $n\in\mathbb{N}$ such that $H^{(n)}=\emptyset$ or $\forall n\in\mathbb{N}\ H^{(n)}\ne\emptyset$.
We can define a mean in the first case. 

\begin{df}Let $H\subset\mathbb{R}$. Let $lev(H)=n\in\mathbb{N}$ if $H^{(n+1)}=\emptyset$ and $H^{(n)}\ne\emptyset$. Otherwise let $lev(H)=+\infty$.
\end{df}

\begin{df}Let $H\subset\mathbb{R},\ lev(H)=n$. Set ${\cal{M}}^{acc}(H)=\A\big(H^{(n)}\big)$. 
\end{df}

In this sense we may say that the last level accumulation points determine the mean and nothing else. Roughly speaking the last accumulation points store the only "weights" of the set.

\begin{lem}\label{llevcup}$lev(H\cup K)=\max\{lev(H),lev(K)\}$.
\end{lem}
\begin{proof}It is known that $(H\cup K)'=H'\cup K'$. From that by induction we get that $(H\cup K)^{(n)}=H^{(n)}\cup K^{(n)}$. Which implies that $(H\cup K)^{(n)}=\emptyset$ iff $H^{(n)}=\emptyset$ and $ K^{(n)}=\emptyset$. 

Now let $m=\max\{lev(H),lev(K)\}$. Then $(H\cup K)^{(m)}\ne\emptyset$ and $(H\cup K)^{(m+1)}=\emptyset$ which gives the statement.
\end{proof}

\begin{lem}\label{llevgr}If $lev(H)<lev(K)$ then ${\cal{M}}^{acc}(H\cup K)={\cal{M}}^{acc}(K)$.
\end{lem}
\begin{proof}By \ref{llevcup} $lev(H\cup K)=lev(K)$. Then $(H\cup K)^{(lev(K))}=K^{(lev(K))}$.
\end{proof}

\begin{lem}$lev(H\cap K)\leq\min\{lev(H),lev(K)\}$.
\end{lem}
\begin{proof}It is known that if $A\subset B$ then $A'\subset B'$ and then by induction $A^{(n)}\subset B^{(n)}$. Apply it for $H\cap K$ and $H$ and then $K$.
\end{proof}

\begin{thm}${\cal{M}}^{acc}$ is strongly internal, finite-independent, shift invariant, symmetric, homogeneous, convex, closed and accumulated generalized mean. 
\end{thm}
\begin{proof}First of all we remark that the definition of ${\cal{M}}^{acc}(H)$ makes sense since $H^{(n)}$ is finite when $H^{(n+1)}=\emptyset$.

${\cal{M}}^{acc}$ is strongly internal because $\varliminf H=\min H',\ \varlimsup H=\max H'$. This gives that ${\cal{M}}^{acc}$ is a generalized mean.

${\cal{M}}^{acc}$ is finite-independent since $H'$ does not change if remove or add finitely many points to $H$. 

It is shift invariant, symmetric, homogeneous since the accumulation operator has the same properties.

To verify convexity let $I$ be a closed interval, ${\cal{M}}^{acc}(H)\in I,\ L\subset I,\ L,L\cup H\in Dom\ {\cal{M}}^{acc}$. Let $lev(H)=n, lev(L)=k$. Now we have three cases: $n<k$, $n>k$, $n=k$. Using \ref{llevgr} the first two are obviously implies that ${\cal{M}}^{acc}(H\cup L)\in I$. For the third \ref{laconvex} gives the statement.

${\cal{M}}^{acc}$ is closed because $cl(H)'=H'$. 

${\cal{M}}^{acc}$ is accumulated since $lev(H')=lev(H)-1$ and \[H^{(lev(H))}=(H')^{(lev(H)-1)}=(H')^{(lev(H'))}.\qedhere\]
\end{proof}

\begin{thm}If either $lev(H)\ne lev(K)$ or $lev(H)=lev(K)=n$ and $H^{(n)}\cap K^{(n)}=\emptyset$ then ${\cal{M}}^{acc}(H\cup K)\in[{\cal{M}}^{acc}(H),{\cal{M}}^{acc}(K)]$.
\end{thm}
\begin{proof}The first case is obvious. For the second case apply $$\A(A\cup B)=\frac{a\A(A)+b\A(B)}{a+b}$$ when $A\cap B=\emptyset$ and $a=|A|,b=|B|$.
\end{proof}

\begin{ex}$H\cap K=\emptyset$ does not imply that ${\cal{M}}^{acc}(H\cup K)\in[{\cal{M}}^{acc}(H),{\cal{M}}^{acc}(K)]$.
\end{ex}
\begin{proof}To show that it is easy to construct sets such that $H\cap K=\emptyset$ and $H'=\{-2,-1,3\},\ K'=\{-1,1\}$. Then ${\cal{M}}^{acc}(H)={\cal{M}}^{acc}(K)=0$ while ${\cal{M}}^{acc}(H\cup K)=\A(\{-2,-1,1,3\})=\frac{1}{4}$.
\end{proof}





\subsection{Means by ideals}\label{ssmbi}
Let us recall the definition of an ideal. ${\cal{I}}\subset P(\mathbb{R})$ is an ideal if $A,B\in {\cal{I}}$ implies that $A\cup B\in{\cal{I}}$ and $B\in{\cal{I}},A\subset B$ implies that $A\in{\cal{I}}$.

\begin{df}Let ${\cal{I}}$ be an ideal. We call ${\cal{I}}$

shift invariant if $H\in{\cal{I}},x\in\mathbb{R}$ implies $H+x\in{\cal{I}}$,

symmetric if $H\in{\cal{I}},x\in\mathbb{R}$ implies $\{x+y:x-y\in H\}\in{\cal{I}}$,

homogeneous if $H\in{\cal{I}},\alpha\in\mathbb{R}$ implies $\alpha H\in{\cal{I}}$.
\end{df}

Evidently the regularly used ideals (e.g. finite sets, countable sets, category 1 sets, sets with Lebesgue measure 0) all have these properties.

\begin{df}Let ${\cal{I}}$ be an ideal and $H\subset\mathbb{R},\ H\notin{\cal{I}}$ be bounded. Set $\varlimsup^{\cal{I}} H=\inf\{x: H^{x+}\in{\cal{I}}\}$. Similarly $\varliminf^{\cal{I}} H=\sup\{x:H^{x-}\in{\cal{I}}\}$.
\end{df}

If ${\cal{I}}=\{\emptyset\}$ then $\varlimsup^{\cal{I}}=\sup,\varliminf^{\cal{I}}=\inf$. If ${\cal{I}}=\{$finite sets$\}$ then $\varlimsup^{\cal{I}}=\varlimsup,\varliminf^{\cal{I}}=\varliminf$. If ${\cal{I}}=\{$countable sets$\}$ then $\varliminf^{\cal{I}},\varlimsup^{\cal{I}}$ are the minimal/maximal consendation points of $H$. If ${\cal{I}}=\{$sets with Lebesgue measure 0$\}$ then $\varliminf^{\cal{I}},\varlimsup^{\cal{I}}$ are the inf/sup of Lebesgue density points of $H$.

\begin{prp}If ${\cal{I}}$ is a $\sigma$-ideal, $H\notin{\cal{I}}$ then $\varlimsup^{\cal{I}} H=\min\{x:H^{x+}\in{\cal{I}}\}$, $\varliminf^{\cal{I}} H=\max\{x:H^{x-}\in{\cal{I}}\}$. \qed
\end{prp}

\begin{prp}\label{pidsu}If ${\cal{I}}_1\subset{\cal{I}}_2,\ H\notin{\cal{I}}_2$ then $\varliminf^{{\cal{I}}_1} H\leq\varliminf^{{\cal{I}}_2} H\leq\varlimsup^{{\cal{I}}_2} H\leq \varlimsup^{{\cal{I}}_1} H $. \qed
\end{prp}

\begin{df}If ${\cal{I}}$ is an ideal, $H\notin{\cal{I}}$ then ${\cal{M}}^{\cal{I}}(H)=\frac{\varliminf^{\cal{I}}(H)+\varlimsup^{\cal{I}}(H)}{2}$.
\end{df}

\begin{thm}\label{thid}Let ${\cal{I}}$ is an ideal. Then ${\cal{M}}^{\cal{I}}$ is a monotone, convex generalized mean.
If ${\cal{I}}$ is shift invariant, symmetric, homogeneous then the mean ${\cal{M}}^{({\cal{I}})}$ has all these properties as well. If $\{$finite sets$\}\subset{\cal{I}}$ then it is finite-independent and strong internal.
\end{thm}

\begin{proof}$\inf H\leq\varliminf^{\cal{I}}(H)\leq\sup H,\ \inf H\leq\varlimsup^{\cal{I}}(H)\leq\sup H$ gives that ${\cal{M}}^{\cal{I}}$ is internal i.e. a mean.

If $\sup H_1\leq\inf H_2$ then $\varliminf^{\cal{I}}(H_1)\leq\varliminf^{\cal{I}}(H_1\cup H_2),\ \varlimsup^{\cal{I}}(H_1)\leq\varlimsup^{\cal{I}}(H_1\cup H_2)$ which yields that ${\cal{M}}^{\cal{I}}(H_1)\leq{\cal{M}}^{\cal{I}}(H_1\cup H_2)$. The other part of monotonicity can be handled similarly.

To verify convexity let $I$ be a closed interval, ${\cal{M}}^{\cal{I}}(H)\in I,\ L\subset I,\ L,L\cup H\in Dom\ {\cal{M}}^{\cal{I}}$. Clearly if $x>\max I$ then $H^{x+}\in\cal{I}$ and because of $L\subset I$ we get that $(H\cup L)^{x+}\in\cal{I}$ which gives that $\varlimsup^{\cal{I}}(H\cup L)\leq\max I$. The other inequality is similar.

If ${\cal{I}}$ is shift invariant, symmetric, homogeneous then so are $\varliminf^{\cal{I}},\varlimsup^{\cal{I}}$ and then so is ${\cal{M}}^{\cal{I}}$.

If $\{$finite sets$\}\subset{\cal{I}}$ then evidently ${\cal{M}}^{\cal{I}}$ is finite-independent hence it is strong internal by \ref{pifiisi}.
\end{proof}

If ${\cal{I}}=\{$finite sets$\}$ then ${\cal{M}}^{({\cal{I}})}={\cal{M}}^{lis}$. If ${\cal{I}}=\{\emptyset\}$ then ${\cal{M}}^{({\cal{I}})}(H)=\frac{\inf H+\sup H}{2}$ that is clearly not strong internal.

\begin{df}Let $({\cal{I}}_n)$ be a sequence of ideals such that ${{\cal{I}}_0}=\{$finite sets$\}$ and ${{\cal{I}}_0}\subset{{\cal{I}}_1}\subset{{\cal{I}}_2}\subset\dots$. The mean associated to this sequence is definied by
$${\cal{M}}^{({\cal{I}}_n)}(H)=\begin{cases}
\A(H) & \text{ if } H \text{ is finite }\\
\frac{\varliminf^{{\cal{I}}_n} H+\varlimsup^{{\cal{I}}_n} H}{2} & \text{ if } H\in{{\cal{I}}_{n+1}}-{{\cal{I}}_n}\\
\lim_{n\to\infty}\frac{\varliminf^{{\cal{I}}_n} H+\varlimsup^{{\cal{I}}_n} H}{2} & \text{ if } H\notin\bigcup^{\infty}_0{{\cal{I}}_n}.\\
\end{cases}$$
\end{df}

\begin{rem}(1) Because of Proposition \ref{pidsu} the limit in the last defining line always exists. 

(2) The definition works for a finite sequence of ideals as well (simply set ${{\cal{I}}_n}={{\cal{I}}_k}$ if $n\geq k$ for a certain $k$).

(3) We can omit the condition that  ${{\cal{I}}_0}=\{$finite sets$\}$. In that case ${\cal{M}}^{({\cal{I}}_n)}$ remains undefined for infinite sets $H\in{{\cal{I}}_1}$.
\end{rem}

The next theorem can be proved like \ref{thid}.

\begin{thm}Let $({\cal{I}}_n)$ be a sequence of ideals such that ${{\cal{I}}_0}=\{$finite sets$\}$ and ${{\cal{I}}_0}\subset{{\cal{I}}_1}\subset{{\cal{I}}_2}\subset\dots$. Then ${\cal{M}}^{\cal{I}}$ is a monotone generalized mean. If ${\cal{I}}_n$ is shift invariant, symmetric, homogeneous then the mean ${\cal{M}}^{({\cal{I}}_n)}$ has all these properties as well.\qed
\end{thm}

\section{Properties of generalized $Avg$}
We can generalize $Avg$ in the following way.

\begin{df}Let $\mu^s$ denote the s-dimensional Hausdorff measure ($0\leq s\leq 1$). If $0<\mu^s(H)<+\infty$ and $H$ is $\mu^s$ measurable (i.e. $H$ is an s-set) then set $$Avg(H)=\frac{\int\limits_H x\ d\mu^s}{\mu^s(H)}.$$
If for a given s we restrict $Avg$ for s-sets then we will use the notation $Avg^s$. 
\end{df}

Clearly $Avg=Avg^0=\A$ for finite sets and we get back the original definition of $Avg^1$ for sets with positive Lebesgue measure (Definition \ref{davg}).

\begin{thm}\label{tavgsish}$Avg$ is shift invariant, symmetric, homogeneous.
\end{thm}
\begin{proof}All properties are a consequence of the theorem on integral by substitution. Let us see them one by one.

We show that $Avg$ is shift invariant: Let $h(x)=x+r\ (r\in\mathbb{R}),\ H$ be an s-set ($0\leq s\leq 1$). Then $$\frac{\int\limits_{h(H)}x d\mu^s}{\mu^s(h(H))}=\frac{\int\limits_{H}x\circ h(x) d\mu^s}{\mu^s(H)}=\frac{\int\limits_{H}x+r d\mu^s}{\mu^s(H)}=\frac{\int\limits_{H}x d\mu^s}{\mu^s(H)}+\frac{\int\limits_{H}r d\mu^s}{\mu^s(H)}=Avg(H)+r$$
where we also used the straightforward fact that $\mu^s(H+r)=\mu^s(H)$.

We prove that $Avg$ is symmetric: By shift invariance it is enough to handle the case when $H$ is symmetric for 0: if $H$ is symmetric for $p\in\mathbb{R}$ then $Avg(H-p)+p=Avg(H)$ and $H-p$ is symmetric for 0 and if $Avg(H-p)=0$ then it would give that $Avg(H)=p$. Let $h(x)=-x$. Then $$\int\limits_{H^{0-}}x d\mu^s=\int\limits_{h^{-1}(H^{0-})}x\circ h(x) d\mu^s=\int\limits_{H^{0+}}-x d\mu^s=-\int\limits_{H^{0+}}x d\mu^s$$  which implies that $Avg(H)=0$.

Now we verify that $Avg$ is homogeneous: Let $h(x)=\alpha x\ (\alpha\in\mathbb{R})$. Then $$\frac{\int\limits_{h(H)}x d\mu^s}{\mu^s(h(H))}=\frac{\int\limits_{h(H)}x d\mu^s}{\int\limits_{h(H)}1 d\mu^s}=\frac{\int\limits_{H}x\circ h(x)\cdot\alpha d\mu^s}{\int\limits_{H}1\circ h(x)\cdot\alpha d\mu^s}=\frac{\alpha^2\int\limits_{H}x d\mu^s}{\alpha\int\limits_{H}1 d\mu^s}=\alpha Avg(H).\qedhere$$
\end{proof}

We can show now that $Avg$ is a generalized mean. For that we prove that it is strong-internal in the following stronger sense.

\begin{prp}Let $H\subset\mathbb{R}$ be a bounded s-set ($0\leq s\leq 1$). Then $\varliminf^{\cal{I}} H<Avg(H)<\varlimsup^{\cal{I}} H$ where ${\cal{I}}=\{H\subset\mathbb{R}:\mu^s(H)=0\}$.  
\end{prp}
\begin{proof}By being shift invariant we can assume that $\varliminf^{\cal{I}} H=0$. We have to prove that $Avg(H)>0$ that is equivalent with $\int\limits_H x\ d\mu^s>0$.
Clearly there is $n\in\mathbb{N}$ such that $\mu^s(H^{\frac{1}{n}+})>0$. Then $0<\frac{1}{n}\mu^s(H^{\frac{1}{n}+})\leq\int\limits_{H^{\frac{1}{n}+}} x\ d\mu^s\leq\int\limits_H x\ d\mu^s.$

The other inequality can be handled similarly.
\end{proof}

\begin{lem}\label{lavgdsu}If $H_1,H_2$ are s-sets and $H_1\cap H_2=\emptyset$ then $$Avg(H_1\cup H_2)=\frac{\mu^s(H_1)Avg(H_1)+\mu^s(H_2)Avg(H_2)}{\mu^s(H_1)+\mu^s(H_2)}$$
\end{lem}
\begin{proof}$$Avg(H_1\cup H_2)=\frac{\int\limits_{H_1}x d\mu^s+\int\limits_{H_2}x d\mu^s}{\mu^s(H_1)+\mu^s(H_2)}=\frac{\mu^s(H_1)Avg(H_1)+\mu^s(H_2)Avg(H_2)}{\mu^s(H_1)+\mu^s(H_2)}.\qedhere$$
\end{proof}

\begin{thm}$Avg$ is strong monotone for s-sets with $s>0$.
\end{thm}
\begin{proof}Let $H_1$ be an $s_1$-set, $H_2$ be an $s_2$-set ($0<s_1,s_2\leq 1$) and $p=\varlimsup H_1\leq\varliminf H_2=r$. Evidently $\mu^{s_1}(H_1^{p+})=0$ and $\mu^{s_2}(H_2^{r-})=0$. Hence $Avg(H_1)\leq p$ and $Avg(H_2)\geq r$.

If $s_1<s_2$ then $Avg(H_1\cup H_2)=Avg(H_2)=r\geq p\geq Avg(H_1)$. 

If $s_2<s_1$ then $Avg(H_1\cup H_2)=Avg(H_1)=p\leq r\leq Avg(H_2)$. 

If $s_1=s_2=s$ then we can assume that $H_1,H_2$ are disjoint because removing a set with 0 measure does not change $Avg$ i.e. $Avg(H_1)=Avg(H_1-H_1^{p+})$ and $Avg(H_2)=Avg(H_2-H_2^{r-})$ . Then by \ref{lavgdsu} $$Avg(H_1\cup H_2)=\frac{\mu^s(H_1)Avg(H_1)+\mu^s(H_2)Avg(H_2)}{\mu^s(H_1)+\mu^s(H_2)}$$
which implies that $Avg(H_1)\leq Avg(H_1\cup H_2)\leq Avg(H_2)$.
\end{proof}

\begin{ex}$Avg$ is not closed and not accumulated either. 
\end{ex}
\begin{proof}Let $H=[0,1]\cup([1,2]\cap\mathbb{Q})$. Then $Avg(H)=Avg^1(H)=0.5$ while $Avg(cl(H))=Avg^1(cl(H))=1$.
\end{proof}

\begin{ex}Symmetry gives $Avg(C)=\frac{1}{2}$ where $C$ is the Cantor set.
\end{ex}

\begin{thm}\label{lint}$Avg$ is convex.  
\end{thm}
\begin{proof}Let $I$ be a closed interval, $Avg(A)\in I,\ C\subset I,\ C,C\cup A\in Dom\ {Avg}$. Let A be an s-set, C be an r-set ($0\leq s,r\leq 1$). 

If $s<r$ then $A\cup C$ is an r-set and $Avg(A\cup C)=Avg^r(A\cup C)=Avg^r(C)=Avg(C)\in I$. If $r<s$ then $A\cup C$ is an s-set and $Avg(A\cup C)=Avg^s(A\cup C)=Avg^s(A)=Avg(A)\in I$.

Let now $s=r$. If $\mu^s(C-A)=0$ then the statement is obvious. Let us suppose $\mu^s(C-A)>0$. By \ref{lavgdsu}
$$Avg(A\cup C)=Avg(A\cup^* (C-A))=$$
$$Avg(A)\frac{\mu^s(A)}{\mu^s(A)+\mu^s(C-A)}+Avg(C-A)\frac{\mu^s(C-A)}{\mu^s(A)+\mu^s(C-A)}\in I$$ 
because it is a convex combination of $Avg(A)$ and $Avg(C-A)$ and both are in $I$.
\end{proof}

\section{Mean by $\epsilon$-neighbourhoods of the set}
We are going to approximate the set by $\epsilon$-neighbourhoods and as they have positive Lebesgue measure, take $Avg$ of those as an approximation of the mean of the set.

Let us recall some usual notation first.

If $H\subset\mathbb{R},\ \epsilon>0$ we use the notation $S(H,\epsilon)=\bigcup_{x\in H}S(x,\epsilon)$ where $S(x,\epsilon)=\{y:|x-y|<\epsilon\}$. Clearly $S(H,\epsilon)=\{y:\exists x\in H\ |x-y|<\epsilon\}$.

\begin{df}Let $H\subset\mathbb{R}$ arbitrary. Set
$$LAvg(H)=\lim_{\delta\to 0+0}Avg(S(H,\delta))$$ if the limit exists. 
\end{df}

\begin{prp}\label{pcll}$LAvg(H)=LAvg(cl(H))$ i.e. $LAvg$ is closed. 
\end{prp}
\begin{proof}It follows from the fact that $S(cl(H),\delta)=S(H,\delta)$.
\end{proof}

This shows that $Avg(H)\ne LAvg(H)$ in general since $Avg$ is not closed.

\begin{thm}\label{tfslavg}Let $H\subset\mathbb{R}$ be a finite set. Then $LAvg(H)=\A(H)$.
\end{thm}
\begin{proof}Let $\delta<\frac{1}{2}\min\{|x-y|:x,y\in H,x\ne y\}$. Then $$AvgS(H,\delta)=\frac{\sum_{x_i\in H}2\delta x_i}{|H|2\delta}=\frac{\sum_{x_i\in H} x_i}{|H|}=\A(H).\eqno\qedhere$$
\end{proof}

\begin{thm}\label{plavgfi}$LAvg$ is finite-independent for infinite sets. 
\end{thm}
\begin{proof}It is enough to prove that for a single point $p$ since from that by induction we can get the statement. If $p$ is an accumulation point then we are done by \ref{pcll}. Let us assume that $p$ is an isolated point and $S(p,\delta_0)\cap S(H-\{p\},\delta_0)=\emptyset$. $H$ is infinite and bounded hence contains an infinite sequence $(h_n)\subset H$ consisting of distinct points. It is enough to show that $\lim\limits_{\delta\to 0+0}\frac{\lambda(S(p,\delta))}{\lambda( S(H,\delta))}=0$ because from that the statement follows since 
$$Avg( S(H,\delta))=\frac{\lambda( S(p,\delta))}{\lambda( S(H,\delta))}Avg( S(p,\delta))+\frac{\lambda(S(H-\{p\},\delta))}{\lambda( S(H,\delta))}Avg( S(H-\{p\},\delta))$$
whenever $\delta<\delta_0$ i.e. when $S(p,\delta)\cap S(H-\{p\},\delta)=\emptyset$.

For that let $K>0$. Find $L\subset H$ such that $|L|=K$. Then find $\delta_1<\delta_0$ such that $l_1,l_2\in L$ implies that $|l_1-l_2|<2\delta_1$. Let $\delta<\delta_1$. Then 
$$\frac{\lambda(S(p,\delta))}{\lambda( S(H,\delta))}<\frac{2\delta}{\lambda( S(L,\delta))}=\frac{2\delta}{2\delta K}=\frac{1}{K}$$
which gives that we required when $K\to\infty$.
\end{proof}

\begin{prp}\label{plavgsi}$LAvg$ is strongly internal. 
\end{prp}
\begin{proof}By \ref{pifiisi} and \ref{plavgfi} it is enough to show that $LAvg$ is internal. 

For that let $m=\inf H,\ \epsilon>0$. If $\delta<\epsilon$ then $Avg(S(H,\delta))>m-\epsilon$. Because it is true for all $\epsilon$ then $LAvg(H)\geq m$. The other inequality is similar.
\end{proof}

\begin{thm}$LAvg$ is shift invariant, symmetric and homogeneous. 
\end{thm}
\begin{proof}Shift invariance comes from $S(H+r,\delta)=S(H,\delta)+r$ and $Avg$ being shift invariant (\ref{tavgsish}). 

Symmetry follows from $S(-H,\delta)=-S(H,\delta)$ and $Avg$ being symmetric (\ref{tavgsish}). 

For proving that $LAvg$ is homogeneous let $\alpha\in\mathbb{R}$. Then $$Avg(S(\alpha H,\delta))=Avg(\alpha S(H,\frac{1}{\alpha}\delta))=\alpha Avg(S(H,\frac{1}{\alpha}\delta)).$$ When $\delta\to 0+0$ then the left hand side tends to $LAvg(\alpha H)$ while the right hand side tends to $\alpha LAvg(H)$.
\end{proof}


\begin{lem}\label{lc0}Let $H\subset\mathbb{R}$ be compact. Then $\forall\epsilon>0\ \exists\delta_0>0$ such that $\delta<\delta_0$ implies $\lambda(S(H,\delta))<\lambda(H)+\epsilon$.  
\end{lem}
\begin{proof} For $\frac{\epsilon}{2}$ there are open intervals $I_i\ (i\in\mathbb{N})$ such that $H\subset\bigcup_1^{\infty}I_i$ and $\sum_1^{\infty}\lambda(I_i)<\lambda(H)+\frac{\epsilon}{2}$. $H$ being compact yields that finitely many covers $H$ as well, e.g. $H\subset\bigcup_1^n I_i$. If we set $\delta_0=\frac{\epsilon}{4n}$ then $\delta<\delta_0$ implies that $\lambda(S(H,\delta))\leq \sum_1^n\lambda(S(I_i,\delta))<\lambda(H)+\frac{\epsilon}{2}+2n\frac{\epsilon}{4n}=\lambda(H)+\epsilon$.
\end{proof}

\begin{thm}\label{pcala}Let $H\subset\mathbb{R}$ be bounded, Lebesgue measurable and $\lambda(H)>0$. Then $Avg(H)=LAvg(H)$ iff $\lambda(cl(H)-H)=0$ or $Avg(cl(H)-H)=Avg(H)$.
\end{thm}
\begin{proof}Let us assume first that $\lambda(cl(H)-H)=0$ i.e. $\lambda(cl(H))=\lambda(H)$. Then clearly $Avg(H)=Avg(cl(H))$ and by \ref{pcll} $LAvg(cl(H))=LAvg(H)$. Hence it is enough to prove the statement for compact sets.

Let $\epsilon>0$ be given. By \ref{lc0} $\forall\epsilon_0>0\ \exists\delta_0>0$ such that $\delta_0<1,\delta<\delta_0$ implies $\lambda(H)\leq\lambda(S(H,\delta))<\lambda(H)+\epsilon_0$. Let $K=\sup H+1$. Then

$|Avg(H)-Avg(S(H,\delta))|=\left\lvert\frac{\int\limits_H x d\lambda}{\lambda(H)}-\frac{\int\limits_{S(H,\delta)} x d\lambda}{\lambda(S(H,\delta))}\right\lvert=$
$\left\lvert\frac{\int\limits_H x d\lambda}{\lambda(H)}-\frac{\int\limits_{H} x d\lambda+\int\limits_{S(H,\delta)-H} x d\lambda}{\lambda(S(H,\delta))}\right\lvert\leq$
$\left\lvert\int\limits_H x d\lambda\right\lvert\left\lvert\frac{1}{\lambda(H)}-\frac{1}{\lambda(S(H,\delta))}\right\lvert+\left\lvert\frac{\int\limits_{S(H,\delta)-H} x d\lambda}{\lambda(S(H,\delta))}\right\lvert\leq$
$K\lambda(H)\left\lvert\frac{\lambda(S(H,\delta))-\lambda(H)}{\lambda(H)\lambda(S(H,\delta))}\right\lvert+\left\lvert\frac{\lambda(S(H,\delta)-H)K}{\lambda(S(H,\delta))}\right\lvert\leq$ 
$\frac{2K}{\lambda(H)}|\lambda(S(H,\delta))-\lambda(H)|<\epsilon$ if $\epsilon_0<\frac{\epsilon\lambda(H)}{2K}$. 

Now assume that $\lambda(cl(H)-H)\ne 0,Avg(cl(H)-H)=Avg(H)$. Then 
\begin{equation}\label{eqlavg}Avg(cl(H))=\frac{\lambda(cl(H)-H)Avg(cl(H)-H)+\lambda(H)Avg(H)}{\lambda(cl(H)-H)+\lambda(H)}=Avg(H).
\end{equation}
Assume now that $Avg(H)=LAvg(H)$. Then apply the first assertion for $cl(H)$. We get $Avg(cl(H))=LAvg(cl(H))=LAvg(H)$ which yields $Avg(cl(H))=Avg(H)$. Then (\ref{eqlavg}) gives the statement.
\end{proof}


%

\begin{ex}\label{elavg1}Let $L=\{\frac{1}{n}:n\in\mathbb{N}\}\cup\{2+\frac{1}{2^n}:n\in\mathbb{N}\}$. Then $LAvg(L)=0$.
\end{ex}
\begin{proof}Let $L_1=\{\frac{1}{n}:n\in\mathbb{N}\}, L_2=\{2+\frac{1}{2^n}:n\in\mathbb{N}\}$.

Let $\delta>0$. Let us estimate where the $\delta$ surroundings $S(x,\delta)$ intersect each other on points of $L_1$ and $L_2$. They intersect on points of $L_1$ when $\frac{1}{n-1}-\frac{1}{n}<2\delta$. It is $n-1>\frac{1}{\sqrt{2\delta}}$. They intersect on points of $L_2$ when $\frac{1}{2^{n-1}}-\frac{1}{2^n}<2\delta$. It is $n>-\log_2{2\delta}$. Then
$$\int\limits_{S(L_1,\delta)} x\ d\lambda<(\sqrt{2\delta}+2\delta))\sqrt{\frac{\delta}{2}}+2\delta(1+\frac{1}{2}+\dots+\frac{1}{n-2})<$$
$$<(\sqrt{2\delta}+2\delta))\sqrt{\frac{\delta}{2}}+2\delta(1-\log\sqrt{2\delta})<2\delta(1-\log\sqrt{2\delta})$$
if $\delta$ is small enough.
$$\lambda(S(L_1,\delta))=\sqrt{2\delta}+2\delta+(n-1-1)2\delta>\sqrt{2\delta}+2\delta+\sqrt{2\delta}-2\delta=2\sqrt{2\delta}$$

$$\int\limits_{S(L_2,\delta)} x\ d\lambda<(2\delta+2\delta)(2+\delta)+2\delta(2+\frac{1}{2^1}+\dots+2+\frac{1}{2^{n-1}})<$$
$$<(4\delta)(2+\delta)+2\delta(2(n-1)+1)<(4\delta)(2+\delta)+2\delta(1-2\log_2{2\delta})<15\delta(1-\log\sqrt{2\delta})$$
if $\delta$ is small enough.
$$\lambda(S(L_2,\delta))=4\delta+(n-1)2\delta>4\delta+2\delta(-\log_2{2\delta}-1)>2\delta+2\delta(-\log{2\delta})$$

$$0<\frac{\int\limits_{S(L_1,\delta)} x\ d\lambda+\int\limits_{S(L_2,\delta)} x\ d\lambda}{\lambda(S(L_1,\delta))+\lambda(S(L_2,\delta))}<\frac{17\delta(1-\log\sqrt{2\delta})}{2\sqrt{2\delta}+2\delta+2\delta(-\log{2\delta})}=$$
$$=\frac{17\sqrt{\delta}(1-\log\sqrt{2\delta})}{2\sqrt{2}+2\sqrt{\delta}+2\sqrt{\delta}(-\log{2\delta})}\to 0\ \text{if }\delta\to 0+0$$
using that $\lim_{x\to 0+0}x\log{x}=0$.
\end{proof}

\begin{prp}$LAvg$ is not accumulated.
\end{prp}
\begin{proof}The example in \ref{elavg1} shows that since $LAvg(L')=1$ by \ref{tfslavg}.
\end{proof}

\begin{prb}Prove or disprove the conjecture that $LAvg$ is an extension of ${\cal{M}}^{iso}$.
\end{prb}

\section{Mean by evenly distributed sample}

Now we define a mean in a way that we take finite sample points from the set and calculate their arithmetic mean and we consider it as an approximation for the mean. It is important that the sample has to be evenly distributed.

\begin{df}\label{d7}Let $H\subset\mathbb{R},\ a=\inf H,b=\sup H$. If $n\in\mathbb{N},0\leq i\leq n-1$ then set $H_{n,i}=H\cap[a+\frac{i}{n}(b-a),a+\frac{i+1}{n}(b-a))$. Let $\ I_n=\{0\leq i\leq n-1:H_{n,i}\ne\emptyset\}$.

We say that the mean of $H$ is $k={\cal{M}}^{eds}(H)$ if $\forall\epsilon>0\ \exists N\in\mathbb{N}$ such that $n>N, \xi_i\in H_{n,i}\ (i\in I_n)$ implies that $|\A(\{\xi_i:i\in I_n\})-k|<\epsilon$.
\end{df}

\begin{thm}If $H\subset\mathbb{R}$ the following statements are equivalent:

(1) ${\cal{M}}^{eds}(H)=k$

(2) $\forall n\in\mathbb{N}$ we select arbitrary points $\xi_{n,i}\in H_{n,i}\ (i\in I_n)$ then $$\lim_{n\to\infty}\A(\{\xi_{n,i}:i\in I_n\})=k$$

(3) $\forall n\in\mathbb{N}$ we select arbitrary points $\xi_{n,i}\in [a+\frac{i}{n}(b-a),a+\frac{i+1}{n}(b-a))\ (i\in I_n)$ then $\lim_{n\to\infty}\A(\{\xi_{n,i}:i\in I_n\})=k$

(4) $\lim_{n\to\infty}\A(\{a+\frac{i}{n}(b-a):i\in I_n\})=k.$
\end{thm}
\begin{proof}(1)$\Leftrightarrow$(2), (3)$\Rightarrow$(4) are obvious.
Proving (2)$\Leftrightarrow$(4)$\Leftrightarrow$(3) at the same time observe that $|A(\{a+\frac{i}{n}(b-a):i\in I_n\}-A(\{\xi_i:i\in I_n\})|\leq\frac{1}{n}$.
\end{proof}

The following theorem verifies that ${\cal{M}}^{eds}$ is a mean.

\begin{thm}${\cal{M}}^{eds}$ is strongly internal.
\end{thm}
\begin{proof}Let $\epsilon>0$. Then $H\cap(-\infty,\varliminf H-\epsilon)$ is finite. $H$ being infinte implies that $\lim_{n\to\infty}|\{H_{n,i}:H_{n,i}\ne\emptyset,H_{n,i}\subset(\varliminf H-\epsilon,+\infty)\}|=\infty$. This gives ${\cal{M}}^{eds}(H)\geq \varliminf H-\epsilon$ by Lemma \ref{ledsm}. Since $\epsilon$ was arbitrary ${\cal{M}}^{eds}(H)\geq \varliminf H$. Similar argument works for $\varlimsup$.
\end{proof}

\begin{prp}If $H$ is finite then ${\cal{M}}^{eds}(H)=\A(H)$.
\end{prp}
\begin{proof}If $n$ is big enough then each interval contains only one point.
\end{proof}

\begin{thm}${\cal{M}}^{eds}$ is shift-invariant, symmetric and homogeneous.
\end{thm}
\begin{proof}Shift-invariance follows from that $\inf,\sup$ and $\A$ are shift-invariant. 

For symmetry it is enough to show that if $0\leq\inf H$ then ${\cal{M}}^{eds}(-H)=-{\cal{M}}^{eds}(H)$. For that we can choose corresponding points in the way that $\xi'_{n,i}=-\xi_{n,i}$ and then we can refer to that $\A$ is symmetric.

For being homogeneous let $f(x)=\alpha x\ (\alpha\in\mathbb{R}^+)$. Note that $f$ takes partition of $[a,b]$ into partition of $[\alpha a,\alpha b]$ and also it takes associated points into associated points of the other partition. Similarly for $f^{-1}$. For completing the proof we need also that \A\ is homogeneous.
\end{proof}

\begin{prp}${\cal{M}}^{eds}$ is monotone and convex.
\end{prp}
\begin{proof}Both statement is a straightforward consequence of \A\ being monotone and convex.
\end{proof}


\begin{prp}If $H=H_1\cup^* H_2$ where $\lambda(cl(H_1))>0$ and $\lambda(H_2)=0, H_2$ is compact then ${\cal{M}}^{eds}(H)={\cal{M}}^{eds}(H_1)$.
\end{prp}
\begin{proof}Let $I^j_n=\{0\leq i\leq n-1:H_{n,i}\ne\emptyset,H_{n,i}\subset H_j\}\ (j\in\{1,2\})$. Then by \ref{lc0} $\lim_{n\to\infty}\frac{1}{n}|I^2_n|=0$. While $\inf\{\frac{1}{n}|I^1_n|:n\in\mathbb{N}\}>0$ which gives the statement using Lemma \ref{ledsm}.
\end{proof}

The next example shows that we cannot omit compactness.

\begin{ex}${\cal{M}}^{eds}([0,1]\bigcup (\mathbb{Q}\cap[1,2]))={\cal{M}}^{eds}([0,2])=1$ hence $Avg\ne{\cal{M}}^{eds}$.
\end{ex}

\begin{ex}Let $L=\{\frac{1}{k}:k\in\mathbb{N}\}\cup\{2+\frac{1}{2^k}:k\in\mathbb{N}\}$. Then ${\cal{M}}^{eds}(L)=0$.
\end{ex}
\begin{proof}Let $a=0,b=2.5,n\in\mathbb{N}$. Let us estimate $|\{i\in I_n:\frac{i+1}{n}(b-a)\leq 1\}|$ i.e. at least how many points $\xi_{n,i}$ we get that are smaller than 1. We want a lower bound. We can get that if $\frac{1}{k-1}-\frac{1}{k}>\frac{1}{n}$. For that it is sufficient that $k<\sqrt{n}$ hence there are at least $\sqrt{n}$ such points.

Now  let us estimate $|\{i\in I_n:\frac{i}{n}(b-a)\geq 2\}|$ i.e. how many points $\xi_{n,i}$ we get that are greater than 2. We want an upper bound. We can get that if $\frac{1}{2^{k-1}}-\frac{1}{2^k}=\frac{1}{2^k}>\frac{1}{n}$ that is $k<\log_2{n}$.

Now $\lim\limits_{n\to\infty}\frac{\log_2{n}}{\sqrt{n}}=0$ completes the proof by Lemma \ref{ledsm}.
\end{proof}

\begin{thm}${\cal{M}}^{eds}(H)\ne{\cal{M}}^{iso}(H)$.
\end{thm}
\begin{proof}Let $H=\{\frac{1}{2^k}:k\in\mathbb{N}\}\cup\{2+\frac{1}{2^k};2+\frac{1}{2^k}+\frac{1}{2^{2^k}}:k\in\mathbb{N}\}$.

Clearly ${\cal{M}}^{iso}(H)=\frac{0+2+2}{3}=\frac{4}{3}$.

Let us calculate ${\cal{M}}^{eds}(H)$. If we divide $[\inf H,\sup H]$ into $2^n$ subintervals  then what is required in order to see points $2+\frac{1}{2^k},2+\frac{1}{2^k}+\frac{1}{2^{2^k}}$  in separate intervals? It is $\frac{1}{2^{2^k}}>\frac{1}{2^n}$ that is $k<\log_2{n}$. Therefore we get $n+1$ points smaller than 1 (converging to 0) and at most $n+1+\log_2{n}$ points greater than 2 (converging to 2). This gives that ${\cal{M}}^{eds}(H)=1$ by Lemma \ref{ledsm}.
\end{proof}

\smallskip
Similar example could show that ${\cal{M}}^{eds}\not< {\cal{M}}^{iso}$ in general.

\begin{prb}Provide example that shows that ${\cal{M}}^{eds}\ne LAvg$. 
\end{prb}


{\footnotesize

\noindent
Dennis G\'abor College, Hungary 1119 Budapest Fej\'er Lip\'ot u. 70.

\noindent E-mail: losonczi@gdf.hu, alosonczi1@gmail.com\\
}
\end{document}